# Using GRASP Approach and Path Relinking to Minimize Total Number of Tardy Jobs on a Single Batch Processing Machine

Panteha Alipour[1], Purushothaman Damodaran*, Christine Nguyen[2]

[1]Department of Industrial Engineering
Purdue University
West Lafayette, IN 47907, USA
palipour@purdue.edu

*Corresponding author: (815) 753-5660; pdamodaran@niu.edu
Department of Industrial and Systems Engineering
Northern Illinois University
DeKalb, IL 60115, USA

[2]Department of Industrial and Systems Engineering
Northern Illinois University
DeKalb, IL 60115, USA
cnguyen@niu.edu

**Author Bio-Sketch**

Panteha Alipour is currently a PhD student in Industrial Engineering at Purdue University. Her research interest is in the area of optimization and scheduling.

Dr. Purush Damodaran is a Professor and Chair of the Department of Industrial and Systems Engineering at Northern Illinois University. His research interest is in the area of large-scale optimization, scheduling, logistics, and simulation. He has completed many industry funded projects in the areas of optimization, scheduling, logistics, simulation, lean, six sigma, and facilities planning. His research group publishes articles regularly on batch processing machine scheduling problems.

Dr. Christine Nguyen is an Assistant Professor in the Department of Industrial and Systems Engineering at Northern Illinois University. She received her doctorate in Industrial and Systems Engineering from the University of Southern California. Her research interests are in the area of optimization, simulation, logistics, supply chain, and data analytics.

**Abstract:** This paper considers the problem of scheduling a single batch processing machine such that the total number of tardy jobs is minimized. The machine can simultaneously process several jobs as a batch as long as the machine capacity is not violated. The batch processing time is equal to the largest processing time among those jobs in the batch. Two decisions are made to schedule jobs on the batch processing machine, namely grouping jobs to form batches and sequencing the batches on the machines. Both the decisions are interdependent as the composition of the batch affects the processing time of the batch. The problem under study is NP-hard. Consequently, solving a mathematical formulation to find an optimal solution is computationally intensive. A Greedy Randomized Adaptive Search Procedure (GRASP) is proposed to solve the problem under study with the assumption of arbitrary job sizes, arbitrary processing times and arbitrary due dates. A novel construction phase for the GRASP approach is proposed to improve the solution quality. In addition, a path relinking procedure is proposed for solving large-sized problems effectively. The performance of the proposed GRASP approach is evaluated by comparing its results to a commercial solver (which was used to solve the mathematical model) and a construction heuristic. Experimental studies suggest that the solution obtained from the GRASP approach is superior compared to the commercial solver and the construction heuristic.

## 1 Introduction

This research is motivated by a practical application where Printed Circuit Boards (PCBs) are assembled and tested. In electronics manufacturing, after assembly, groups of PCBs are tested in an Environmental Stress Screening (ESS) chamber to detect component and other board-level failures before they are used in the field. ESS is the process of applying environmental stresses in conjunction with functional testing in order to stimulate the failure mechanisms of defects to the point of detection, which could not be detected by visual inspection or in-circuit testing (Zinouri, Muskeyvalley, Damodaran, and Ghrayeb, 2012). The PCBs are typically subjected to vibration, humidity, and thermal cycling or power cycling in an ESS chamber. The ESS chamber is referred to as a Batch Processing Machine (BPM), and the PCBs are represented as jobs in this research to remain consistent with the scheduling literature.

The jobs are assumed to vary in terms of size (i.e. smaller PCBs used in consumer products such as cell phones and larger PCBs in industrial applications such as servers). The processing time, which is the testing time, also varies depending upon the customer's testing requirements. Testing times are prescribed by customers. While testing, the PCBs should be tested at least for the time prescribed by the customer. Customers do not pay for the testing time above the prescribed minimum. Once several jobs are grouped as a batch, the processing time of the batch is equal to the longest processing time of all the jobs in the batch. Two tasks are carried out to schedule a set of jobs on a BPM. The first task is to form batches such that the machine capacity is not violated, and the second task is to determine a schedule that minimizes the total number of tardy jobs. The problem of minimizing the total number of tardy jobs on a single batch processing machine is NP-hard in strong sense (Brucker, Kovalyov, Shafransky, and Werner, 1998), and hence a Greedy Randomized Adaptive Search Procedure (GRASP) with path relinking is proposed in this research. Path relinking implements a novel local search strategy based on different neighbourhood structures defined by path exchanges, and it uses memory based on a local search strategy to find better solutions in less computation time (Aiex, Binato, and Resendel, 2003; Resendel and Ribeiro, 2005).

## 2 Literature review

There are many researchers who studied makespan and completion time objectives while scheduling a BPM. Uzsoy (1994) showed that the single BPM scheduling problem with non-identical job sizes to minimize the total completion time and makespan is NP-hard. He proposed bin-packing-based heuristics for minimizing makespan and has used the Branch-and-Bound (B&B) approach to minimize the completion times. He also developed effective heuristics for minimizing makespan and total completion time. Lee (1999) proposed polynomial and pseudo-polynomial time algorithm to minimize makespan of a single BPM with dynamic job arrivals (i.e. job release times are not equal). The algorithm presented excellent results but had long computational time. Chandru et al. (1993) investigated the minimizing makespan objective in a single BPM by B&B. The set of jobs to be scheduled are grouped into a number of families, where all jobs in the same family have the same processing time. Vélez-Gallego et al. (2011) studied constructive heuristics named as Modified Successive Knapsack Problem (MSKP) to minimize makespan in a single BPM under the assumptions of non-identical job sizes and non-zero job ready times. The heuristic was found to outperform other comparative approaches for test instances with 50 or more jobs.

On-time jobs are the ones which are completed before their due dates, and tardy jobs are the jobs which finish after their due dates. Moore (1968) minimized the total number of tardy jobs on a discrete processing machine $(1| \, |\sum U_j)$ using Hodgson's algorithm. His approach was able to guarantee an optimal solution. Lenstra, Kan and Brucker (1977) showed that scheduling a discrete machine to minimize the total weighted number of tardy jobs with ready time constraint $(1|r_j|\sum w_j U_j)$ is NP-hard. Shanthikumar (1983) considered the one machine scheduling problem with dual objective of minimizing the maximum tardiness with minimum number of tardy jobs. Perez et al. (2005) proposed a heuristic to minimize total weighted tardiness in the diffusion step of semiconductor wafer fabrication process in single BPM with non-identical job sizes.

Batch machine scheduling problems where the batch processing time is given by the processing time of the longest job in the batch have also been addressed. Brucker et al. (1998) present several exact algorithms for minimizing the number of tardy jobs on a batch processing machine $(1|p-batch|\sum U_j)$. They present a dynamic programming approach for minimizing the total number of tardy jobs and show that minimizing total weighted number of tardy jobs is binary NP-hard. They assumed that only a fixed number of jobs are allowed in a batch. They did not consider the job sizes.

Ourari et al. (2009) considered the $1|r_j|\sum U_j$ problem, and they used an original mathematical integer programming formulation where they showed how both good-quality lower and upper bounds can be computed. Sevaux and Dauzère-Pérès (2003) presented the first meta-heuristic (i.e. genetic algorithm) for a discrete processing machine to minimize the total number of tardy jobs with ready time limitation $(1|r_j|\sum w_j U_j)$. M'Hallah and Bulfin (2007) also tried to minimize the weighted number of tardy jobs with release date for single discrete machine $(1|r_j|\sum w_j U_j)$. They gave a formulation to maximize the weighted number of on-time jobs instead of minimizing weighted number of tardy jobs.

Jolai (2005) proved that minimizing the number of tardy jobs with incompatible job families on a single BPM is NP-hard. Also, Liu and Zhang (2008) proved that minimizing the number of tardy jobs on a batch processing machine with incompatible job families is unary NP-hard in strong sense.

Several review papers provide expanded records of different aspects of scheduling. The survey of minimizing total number of late jobs problems is given by Leung (2004). To the best of our knowledge, the problem of minimizing total number of tardy jobs on a batch processing machine with non-identical job sizes, machine capacity, and arbitrary due dates has not been studied so far. A GRASP approach is proposed in this research.

GRASP is a meta-heuristic, which is an iterative process popularized by Feo and Resende (1989). This multi-start method includes two different phases. A classical GRASP starts by producing a feasible solution, which is constructed by a greedy randomized algorithm, and tries to improve it in the local search phase. The procedure is repeated several times; and the local optimum in the neighborhood of the constructed solution is desired. The best overall solution is returned as a result (Feo and Resende, 1989; Feo and Resende, 1995). Resende and Ribeiro (2014) present expanded surveys on GRASP. References are also found that implemented GRASP to solve scheduling problems (Aiex et al., 2003; Armentano and de Araujo, 2006; Damodaran et al., 2013; Damodaran et al., 2011).

## 3 Mathematical formulation

The Mixed Integer Linear Programming (MILP) model for minimizing the total number of tardy jobs on a batch processing machine with non-identical job sizes is presented in this section. The notation used in the formulation is presented below.

Sets

| | |
|---|---|
| $\{j \in J\}$ | Set of jobs |
| $\{b \in B\}$ | Set of batches |

Parameters

| | |
|---|---|
| $p_j$ | Processing time of job *j* |
| $s_j$ | Size of job *j* |
| $d_j$ | Due date of job *j* |
| $S$ | Machine capacity |
| $M$ | A very large positive number |
| $e$ | A very small positive number |

Decision Variables

| | |
|---|---|
| $c_j$ | Completion time of job $j$ |
| $C_b$ | Completion time of $b^{\text{th}}$ batch scheduled |
| $P_b$ | Processing time of $b^{\text{th}}$ batch scheduled |
| $NT_j$ | 1, if job $j$ is tardy; 0, otherwise |
| $X_{jb}$ | 1, if job $j$ is assigned to the $b^{\text{th}}$ batch; 0, otherwise |

The mixed integer linear programming formulation for the problem under study is presented below.

Minimize $\sum_{j \in J} NT_j$ (1)

*Subject to:*

$$\sum_{b \in B} X_{jb} = 1 \qquad \forall j \in J \qquad (2)$$

$$\sum_{j \in J} s_j X_{jb} \leq S \qquad \forall b \in B \qquad (3)$$

$$P_b \geq P_j X_{jb} \qquad \forall j \in J, b \in B \qquad (4)$$

$$C_1 = P_1 \qquad (5)$$

$$C_b \geq C_{b-1} + P_b \qquad \forall b \in B/\{1\} \qquad (6)$$

$$c_j \geq C_b - M(1 - X_{jb}) \qquad \forall b \in B \qquad (7)$$

$$c_j - d_j \leq M(NT_j) \qquad \forall j \in J, b \in B \qquad (8)$$

$$c_j - d_j \geq e - M(1 - NT_j\ ) \qquad \forall j \in J, b \in B \qquad (9)$$

$$C_b, P_b \geq 0 \qquad \forall b \in B \qquad (10)$$

$$c_j \geq 0, NT_j \in \{0,1\} \qquad \forall j \in J \qquad (11)$$

$$X_{jb} \in \{0,1\} \qquad \forall j \in J, b \in B \qquad (12)$$

The objective (1) is to minimize the total number of tardy jobs. Constraint set (2) ensures that each job is assigned to exactly one batch. Constraint set (3) ensures that the total size of all jobs in each batch processed on the machine does not violate the machine capacity. Constraint set (4) determines the processing time of the $b^{th}$ batch processed on the machine. Batch processing time is equal to the longest processing time of all the jobs in the batch. Constraint set (5) ensures the completion time of the first batch is equal to the processing time of the first batch processed in the sequence. Constraint set (6) ensures that completion time of the batch is at least equal to the summation of its processing time and completion time of the previous batch. Constraint set (7) ensures that the completion time of a job is at least equal to the completion time of the batch in which it is processed. Constraint sets (8) and (9) determine whether a job is tardy or not. Constraint sets (10), (11), and (12) impose the non-negativity and binary restrictions on the decision variables. In constraints (7), the parameter $M$ is chosen to be large enough to ensure that the completion time of job $j$ is greater or equal to the completion time of $b^{th}$ batch. For constraints (8) and (9) $M$ and $e$ are chosen to determine whether job $j$ is late or not.

The above formulation can be solved in a commercial mixed integer linear program solver such as IBM ILOG CPLEX. Since the problem is NP hard, as the number of jobs increases, it is expected that the commercial solver will require prohibitively long run time to solve the problem to optimality. The experimental results (see section 6) also indicate that the solver requires a long run time to solve larger problem instances. Consequently, commercial solvers may not be viable to use in a real-life setting. To solve the problem under study more efficiently and effectively, a GRASP approach is proposed in the next section.

## 4 GRASP Solution approach

GRASP is a heuristic iterative sampling technique composed of two phases: a construction and a local search phase. In the construction phase, an iterative algorithm builds a feasible solution by adding one variable (or job) at a time. Each iteration of the algorithm starts with a solution found by means of a randomized greedy heuristic. The solution is later taken as the initial solution of the local search procedure and the procedure is repeated until some criteria of the search are met. The pseudo-code for the GRASP is shown in Figure 1.

---

**BEGIN GRASP**

**for** *j=1,…, Num Times To Run* **do**

*NumberIterNoImprove*← 0

**for** *Iter=1,…, MaxIters* **do**

$x^*$ ← *Greedy Randomized Construction (Seed);*

*x* ← *Local search (Solution);*

**if** *f(x)>f* $(x^*)$ **then**

$x^* \leftarrow x$; $f^* \leftarrow f(x)$; *NumberIterNoImprove* ← 0;

**else**

*NumberIterNoImprove*←*NumberIterNoImprove+1;*

**end if**

**end for**

**end for**

**return ($x^*$);**

**end GRASP**

---

**Figure 1** Pseudo-code for GRASP.

While a greedy heuristic is a construction heuristic that fixes one variable at a time using some deterministic rule, in a randomized greedy heuristic, randomness is introduced to heuristic rules in order to prevent generation of the same solution every time the greedy heuristic is called. The randomized greedy heuristic proposed to find a starting solution in this research is based on several different rules: earliest due date (EDD), shortest processing time, shortest size, largest size, the difference between due date and processing time which is multiplied by size, shortest result of multiplication of size and processing time, shortest result of multiplication of due date and processing time, shortest result of multiplication of due date and size, and random sequence (i.e. random job permutation).

If each job is assigned to its' own unique batch, then the maximum number of batches formed is equal to the number of jobs ( $|B| = |J| = n$ ) considered in an instance. Initially, $n$ batches are opened, and any remaining unused batches are closed at the completion of the algorithm. The jobs are first sequenced based of one of the 10 rules (i.e. shortest processing time, earliest due date, etc. as shown in the pseudo code in Figure 2). The earliest job from the initial sequence generated is assigned to the first available batch (i.e. adding a job should not violate the machine capacity). This procedure is repeated until all jobs are assigned to a batch. After the batches are formed, the batch processing time can be determined. The batch processing time is equal to the largest processing time of all jobs in the batch. The completion time of each job in a batch is equal to the completion time of the batch. The completion time of a batch can be easily computed by summing the processing times of all the batches sequenced (including the processing time of the current batch). The batches are sequenced in the order in which they were created. If the completion time of the job is larger than its due date, then the job would be tardy.

After generating all initial sequences, the sequences are used to find the total number of tardy jobs. In the next step, the pool of path relinking is created, which contains all the initial sequences and their total number of tardy jobs.

Randomization is introduced to the earliest due date, shortest processing time, shortest size and other initial sequences by means of the so-called Restricted Candidate List (RCL). It works as follows: while assigning jobs to the batches, instead of selecting the unscheduled job with the earliest due date first (for EDD sequence), the job is randomly chosen among the first *k* earliest due date jobs. Here *k* is the size of the RCL.

The pseudo-code for the construction phase of GRASP is shown in Figure 2.

---

**BEGIN**

Let $B$ = {Batch$_1$, Batch$_2$,…,Batch$_n$} be the set of $n$ initial empty batches.

Let $J$ be the set of $n$ jobs.

Job parameters – due date ($d_j$), size ($s_j$), and processing time ($p_j$) are given

Compute $A_j = s_j p_j$, $B_j = s_j d_j$, $C_j = s_j(d_j - p_j)$ and $E_j = p_j d_j$.

**For *k* = 1** to **10 do**

**While** $|J| > 0$ **do**

$q_{k=1} \leftarrow argmin_{j\in J}\{d_j\}$ //to create first initial sequence based on EDD sequence

$q_{k=2} \leftarrow argmin_{j\in J}\{s_j\}$ //to create second initial sequence based on smallest sized job first

$q_{k=3} \leftarrow argmax_{j\in J}\{s_j\}$ //to create third initial sequence based on largest sized job first

$q_{k=4} \leftarrow argmin_{j\in J}\{p_j\}$ //to create fourth initial sequence based on smallest processing time job first

$q_{k=5} \leftarrow argmin_{j\in J}\{A_j\}$ //to create fifth initial sequence based on smallest $A_j = p_j.s_j$ first sequence

$q_{k=6} \leftarrow argmin_{j\in J}\{B_j\}$ //to create sixth initial sequence based on smallest $B_j = s_j.d_j$ first sequence

$q_{k=7} \leftarrow argmin_{j\in J}\{C_j\}$ //to create seventh initial sequence based on smallest $C_j = s_j.(d_j\text{-}p_j)$ first sequence

$q_{k=8} \leftarrow argmin_{j\in J}\{E_j\}$ //to create eight initial sequence based on smallest $p_j.d_j$ first sequence

$q_{k=9,k=10} \leftarrow randperm(J)$ //last sequences are random permutations of all jobs

Sequence$_k$ = Sequence$_k \cup \{q_k\}$ //jobs are sequenced in an order based on the value of *k*

**end while**

**end for**

**For *k* = 1** to **10 do**

**While** |Sequence$_k$| > 0 **do**

$q_k \leftarrow$ randomly pick one job from the first $K$ jobs ($K$ is the size of RCL)

**for** $t$ =1 to $n$ **do**

**if** $\sum_{b\in Batch_t} s_b + s_{q_k} \leq S$ **then** // if job $q_k$ can be allocated in batch $t$, then

$Batch_{tk} \leftarrow Batch_{tk} \cup \{q_k\}$ // assign job $q_k$ to batch $t$

$J \leftarrow J \backslash \{q_k\}$ // delete job $q_k$ from the job set $J$

$P_{Batch_{tk}} \leftarrow max_{j\in Batch_{tk}}\{p_j\}$ //batch processing time equal to largest processing job in the batch

**exit for loop**

**end if**

**end while**

$C_{1k} = P_{Batch_{1k}}$ //completion time of first batch is equal to processing time of the first batch

**for** $t$ = 2 to $n$ **do**

$C_{tk} = C_{t-1k} + P_{Batch_{tk}}$

**end for**

Compute the job completion times //job completion time is equal to the completion time of the batch to which the job was assigned

Compute number of tardy jobs //count the number of jobs completed after their due date

**end for**

**Figure 2** Pseudo-code for the construction phase of GRASP.

**Table 1** Data for a 9-job problem instance.

| Jobs # | 1 | 2 | 3 | 4 | 5 | 6 | 7 | 8 | 9 |
|---|---|---|---|---|---|---|---|---|---|
| **Job size** | 17 | 13 | 27 | 7 | 15 | 14 | 27 | 2 | 28 |
| **Job Processing Time** | 19 | 28 | 44 | 14 | 16 | 23 | 37 | 10 | 43 |
| **Job Due Date** | 36 | 35 | 32 | 32 | 34 | 36 | 36 | 37 | 36 |

Table 1 shows the data for a 9-job problem instance used to illustrate the greedy randomized algorithm. The machine capacity is assumed as 40 and the RCL as 3. The jobs are first arranged in EDD order (i.e. list = {3,4,5,2,1,6,7,9,8}). In the classical greedy heuristic, the first unscheduled job from the list is assigned to the first feasible batch. If none of the existing batches can accommodate the selected job, then a new batch is created, and the job is assigned to the new batch. Finally, the job is deleted from the list and the procedure starts again until all the jobs have been assigned to a batch. However, in the randomized greedy heuristic, a job is randomly chosen from the first 3 unscheduled jobs (as RCL = 3) in the list. Table 2 shows how the randomized greedy algorithm is applied to a 9-job problem instance. Instead of choosing job 3 first, the randomized greedy algorithm randomly chooses a job from the top 3 jobs in the list (in this example job 5 was chosen). Since it is the first job chosen, it will fit in batch 1. Second, the job 4 is chosen and it will also fit in batch 1. However, the third job chosen is job 3 which does not fit in batch 1 as its addition would violate the machine capacity. Consequently, job 3 is placed in the second batch. Next, job is 1 chosen from the list and as it fits in batch 1 it is assigned to batch 1. As the total batch size of batch 1 is 39 and there are no jobs with a size of 1, batch 1 can be closed. This procedure is repeated until all the jobs are assigned to a batch. Later the batches are processed (or scheduled) on the machine in the order in which they were formed. Applying the randomized greedy heuristic for this example problem resulted in five batches as shown in Table 2. Batch processing times and batch completion times can be easily determined once the batch composition is known. Figure 2 presents a Gantt chart for the batches formed. For this example, problem, six jobs are tardy (i.e. jobs 3, 2, 7, 8, 9, and 6 are tardy).

**Table 2** Greedy randomized heuristic using EDD job sequence.

| Step | List | Random | Job | Batch size | Batch | Jobs | Batch Processing Time | Batch Completion Time |
|---|---|---|---|---|---|---|---|---|
| **1** | {3,4,5,2,1,6,7,9,8} | 3 | 5 | 15 | 1 | {5,4,1} | 19 | 19 |
| **2** | {3,4,2,1,6,7,9,8} | 2 | 4 | 22 | 2 | {3,2} | 44 | 63 |
| **3** | {3,2,1,6,7,9,8} | 1 | 3 | 27 | 3 | {7,8} | 37 | 100 |
| **4** | {2,1,6,7,9,8} | 2 | 1 | 39 | 4 | {9} | 43 | 143 |
| **5** | {2,6,7,9,8} | 3 | 7 | 27 | 5 | {6} | 23 | 166 |
| **6** | {2,6,9,8} | 1 | 2 | 40 | | | | |
| **7** | {6,9,8} | 3 | 8 | 29 | | | | |
| **8** | {6,9} | 2 | 9 | 28 | | | | |
| **9** | {6} | 1 | 6 | 5 | | | | |

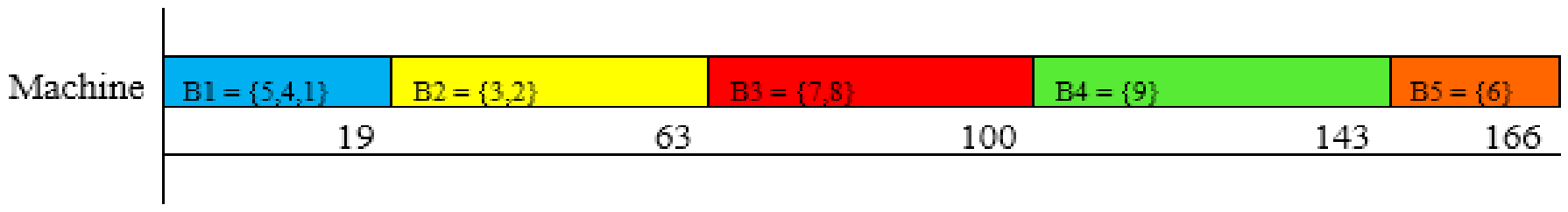

**Figure 2** Gantt chart for greedy randomized algorithm.

The randomized greedy heuristic can be further improved by an observation made in Jolai (2005) and Sevaux and Dauzère-Pérès (2003). At the time of adding a job to the first available batch, if the batch completion time is greater than the due date of the job added, then the job is already tardy. Consequently, adding the job to any batch that follows the first available batch will also make the job tardy. Thus, to minimize the number of tardy jobs, all the batches which contain jobs that are completed before their due dates are processed first on the machine. Table 3 illustrates the improved greedy randomized heuristic for the problem instance shown in Table 1. Figure 3 presents the corresponding Gantt chart. Five jobs are tardy for this example (i.e. jobs 1, 6, 3, 7, and 9).

**Table 3** Improved greedy randomized heuristic for EDD sequence.

| Step | List | Random | Job | Tardy or Non-tardy | Batch size | Batch | Batch | Jobs | Batch processing time | Batch completion time |
|---|---|---|---|---|---|---|---|---|---|---|
| **1** | {3,4,5,2,1,6, 7,9,8} | 2 | 4 | Is not tardy $(c_j < d_j)$ & $(p_j < d_j)$ | 7 | 1 | 1 | {4,2,5, 8} | 28 | 28 |
| **2** | {3,5,2,1,6,7, 9,8} | 3 | 2 | Is not tardy $(c_j < d_j)$ & $(p_j < d_j)$ | 20 | 1 | 2 | {1,6} | 23 | 51 |
| **3** | {3,5 ,1,6,7,9,8} | 2 | 5 | Is not tardy $(c_j < d_j)$ & $(p_j < d_j)$ | 35 | 1 | 3 | {3} | 44 | 95 |
| **4** | {3,1,6,7,9,8} | 2 | 1 | Does not fit in first batch. Adding the job to second batch would make it tardy. $(c_j > d_j)$ | 17 | 2 | 4 | {7} | 37 | 132 |
| **5** | {3,6,7,9,8} | 1 | 3 | Does not fit in batches 1 and 2. This job will be tardy as the processing time is already larger than its due date. $(p_j > d_j)$ | 27 | 3 | 5 | {9} | 43 | 175 |
| **6** | {6,7,9,8} | 2 | 7 | Does not fit in batches 1, 2 and 3. This job will be tardy as the | 27 | 4 | | | | |

| | | | | | | |
|---|---|---|---|---|---|---|
| | | | | processing time is already larger than its due date. $(p_j > d_j)$ | | |
| **7** | {6,9,8} | 3 | 8 | Can fit in first batch. Adding this job to the first batch does not make it tardy. | 37 | 1 |
| **8** | {6,9} | 2 | 9 | Does not fit in batches 1 through 4. This job will be tardy as the processing time is already larger than its due date. $(p_j > d_j)$ | 28 | 5 |
| **9** | {6} | 1 | 6 | Can fit in batch 2. However, it is tardy. $(c_j > d_j)$ | 31 | 2 |

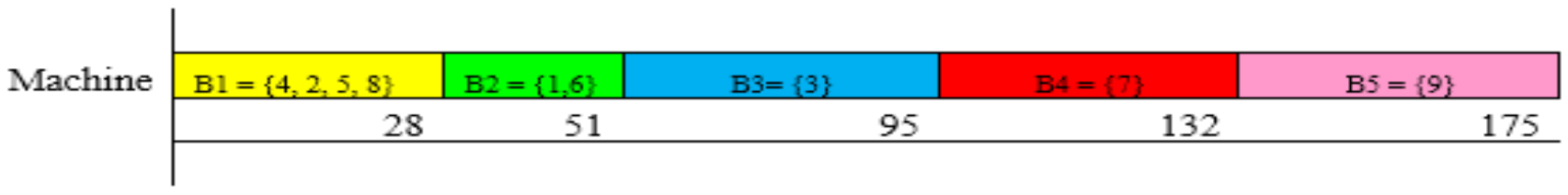

**Figure 4** Gantt chart for improved greedy randomized algorithm.

Three different types of movements define the neighbourhood in the local search procedure: (1) batch interchange, where two batches are selected at random and interchanged. This movement is inspired by Damodaran et al. (2011). Interchanging a batch with longer batch processing time with a batch with shorter processing time located earlier in the sequence may help reduce the completion time of jobs and thereby the number of tardy jobs; (2) job insertion (A), where the job with the largest processing time is marked and discarded from one batch and then inserted into a randomly selected batch or a newly created batch (if the inclusion of the new job would lead to additional tardy jobs); (3) job insertion (B), where jobs with very large processing times compared to the mean of the processing time of all jobs in that batch are removed and inserted to an existing or a new batch. If the processing time of the job is larger than the arbitrary rate of the mean of the processing time of all jobs in that batch ($(p_j > (\alpha + 1))$ * mean of all $p_j$s in that batch ), the job is discarded from that batch and inserted into some other batch or newly created batch. Although this movement is similar to the previous one, this movement determines the job with larger processing time as a candidate for the insertion operation.

Table 4 presents the data for a 9-job instance used to illustrate the various movements described above. The machine capacity is assumed to be 40. The jobs were initially sequenced based on the EDD rule. Four batches were formed. The Gantt chart for the schedule is shown in Figure 5. All the nine jobs are tardy as per this schedule. If batch 2 is interchanged with batch 1, then seven jobs are tardy (i.e. jobs 1, 2, 3, 5, 6, 7, and 8 are tardy).

**Table 4** Data for a 9-job instance.

| **Jobs numbers** | **1** | **2** | **3** | **4** | **5** | **6** | **7** | **8** | **9** |
|---|---|---|---|---|---|---|---|---|---|
| **Job size** | 37 | 18 | 5 | 12 | 9 | 2 | 10 | 4 | 25 |
| **Job Processing Time** | 22 | 4 | 3 | 2 | 24 | 50 | 8 | 5 | 10 |
| **Job Due Date** | 35 | 10 | 12 | 21 | 26 | 15 | 17 | 36 | 24 |

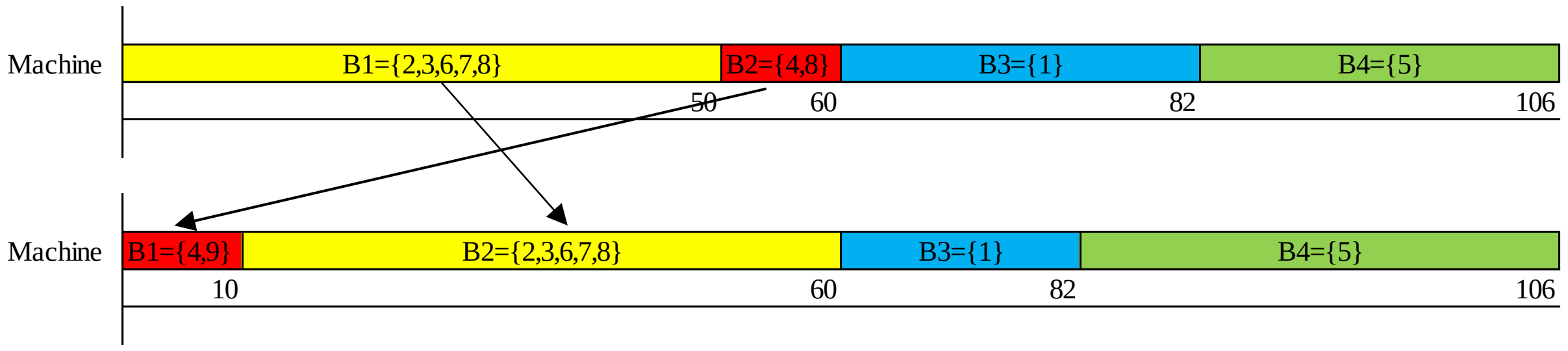

**Figure 5** Batch interchange.

A job insertion type (A) is shown in Figure 6, where job 6 with larger processing time is discarded from batch 1 and inserted into a newly created batch (i.e. batch 5). Here the total number of tardy jobs is 3 compared to the original schedule with 9 tardy jobs. Job 6 is discarded from the first batch and inserted into a newly created batch to avoid making other jobs tardy.

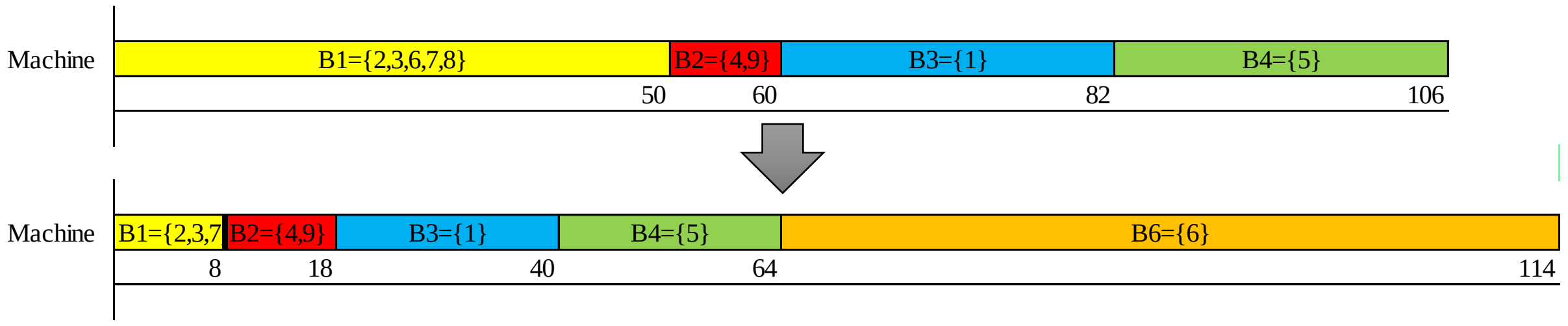

**Figure 6** Job insertion (A).

A job insertion type (B) is shown in Figure 7. This movement is slightly similar to the previous movement; however, here jobs with larger processing time are compared with the mean processing time of the batch. For batch 1, the mean of the processing time of all jobs is $\frac{50+4+3+8+5}{5} = 14$ and the processing time of job 6, which is 50, is larger than the mean of processing time of all jobs (in this example $\alpha$ is assumed to be equal to epsilon). Thus job 6 is removed from the first batch and the batch completion time is reduced (or job completion times are reduced). Again, for the second batch, the mean of processing time of jobs 4 and 9 is 6. Job number 9 with processing time equal to 10 is the job with largest processing time in the second batch (comparing to the mean of processing time of jobs 4 and 9), so it is also removed from the second batch. All in all, jobs 9 and 6 with larger processing times are discarded from batches 1 and 2 and inserted into a newly created batch (batch 5). Here the total number of tardy jobs is 2. Jobs 6 and 9 with largest processing time are removed to avoid making other jobs tardy.

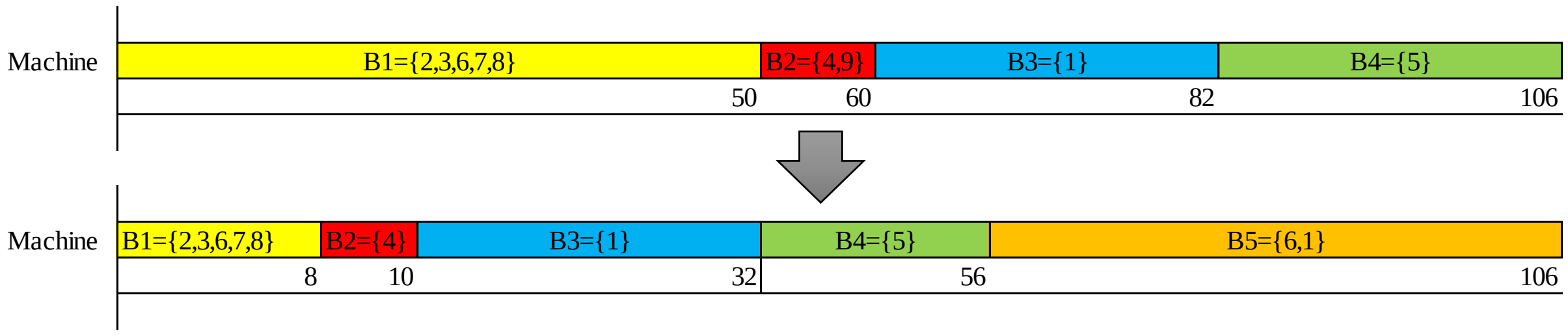

**Figure 7** Job insertion (B).

The approach of path relinking for use within GRASP is developed by Laguna and Marti (1999). A set of high-quality solutions is saved as a pool to provide the initial and guiding solutions for the path relinking procedure. Each iteration of GRASP produces a local optimal solution (a sequence of jobs in the problem under study) and this solution is considered as an initial solution (I). Another good solution, which is the best among other solutions, in the pool is chosen and it is named guiding (G) solution, and a path of solutions linking the initial solution to the guiding solution is constructed by applying a series of changes to the original solution (Resende and Ribeiro, 2005). For instance, let I= (1,3,4,2,5) and G= (5,3,2,4,1). A path relinking I to G is I = (1,3,4,5,2) → (2,3,4,5,1) →(5,3,4,2,1) →(5,3,2,4,1) = G and each of these solutions is evaluated for solution quality (minimum of total number of tardy jobs in the problem under study). The pseudo-code for the proposed path relinking approach is shown in Figure 8.

---

Procedure path relinking

**for** all PR-iteration=1 to max PR-iteration do

R ←Gen Rand Int (# of Elite pool)

Initial Seq←R

Guiding Seq←Best Seq of Elite pool

Set Initial Seq ←first Seq

**while** (Initial Seq ≠ Guiding Seq)

**for** all n= 1 to $n_j$

**if** (Initial Seq (n) ≠ Guiding Seq(n)) then

Find t │ Initial Seq (t) == Guiding Seq(n)

Swap Initial Seq (n) ←initial Seq (t)

**end if**

**end for**

**end while**

**end for**

**END**

---

**Figure 8** Pseudo-code for path relinking.

## 5 Data source and description

In this section, we present the data set that we used to test our approach. In this experiment, different examples of 5 to 100 jobs have been used. The data set was generated using the parameters shown in Table 5. In order to generate job due dates, due date adjustment factor, job ready times, job processing times, and $z_j$ factor were considered. The due date is determined by setting $z_j$ and finding the makespan using a Full Batch LPT (FBLPT) schedule. For creating FBLPT schedule, the jobs are ordered in the non-increasing order of their processing times. Next, the jobs are grouped and placed in batches one by one based on their size and with consideration of capacity restriction.

Throughout this study, small number of jobs data set refers to the data set consisting of 5, 7, 9, 11, 13 and 15 jobs instances, while large number of jobs data set refers to the data set consisting of 50 and 100 jobs instances.

**Table 5** Parameters' description and their values.

| Parameters | Description | Values |
|---|---|---|
| $n$ | Number of jobs (problem size) | 5-100 |
| $m$ | Number of machines | 1 |
| $S_m$ | Machine capacity | 40 |
| $Q$ | Interval for job sizes | [1,30] |
| $R$ | Due date tightness | R=0.5 |
| $T$ | Due date spread out | T=0.3 |
| $r_j$ | Job ready times | Discrete uniform [0,48] |
| $p_j$ | Job processing times | Discrete uniform [8,48] |
| $w_j$ | Job priorities | Discrete uniform [1,11] |
| $d_j$ | Job due dates | $\gamma.(r_j + p_j + z_j)$<br>Where<br>$z_j$=Discrete uniform [$\mu$.(1-R/2), $\mu$.(1+R/2)]<br>$\mu$ = (1-T). $C^*_{max}$<br>$C^*_{max}$=MIN ($r_j$) +FBLPT<br>FBLPT is Full Batch LPT (zero ready times and unit sized jobs) |
| $\Gamma$ | Due date adjustment factor | 0.2, 0.33, 0.5 |

## 6 Experimental study

An experiment was conducted to fine tune GRASP with path relinking parameters (i.e. the size of RCL, number of iterations in GRASP, and number of path relinking iterations). Table 6 presents the different factors and levels considered during preliminary study to tune GRASP parameters. The levels tested were $k$ = {5 for small jobs instances, 10% of number of jobs, 25% of number of jobs and 50% of number of jobs for large jobs instances}, total number of iterations for GRASP= {300, 1000}, and path relinking iterations = {300, 1000}.

**Table 6** Factors and levels.

| Factors | Levels |
|---|---|
| Size of RCL, $k$ | $10\%n$, $25\%n$ and $50\%n$ |
| Number of GRASP iterations | 300, 1000 |
| Number of path relinking iterations | 300, 1000 |

During this process, all 50-job instances were solved using the GRASP algorithm and the commercial solver (i.e. IBM ILOG CPLEX). There were altogether ten 50-job instances assuming a machine capacity $S = 40$; each instance was run for all combinations of the GRASP with path relinking parameters shown in Table 5 (i.e., $10\times 3 \times 2 \times 2 = 120$ experiments). The combination of each parameter which yielded the most improvement in minimizing the total number of tardy jobs is chosen for further experimentation. The average percent improvement in minimizing total number of tardy jobs was computed using Equation (13).

$$\%\ \text{Improvement} = \frac{(\sum \mathrm{U_j}\ \text{from CPLEX} - \sum \mathrm{U_j}\ \text{from GRASP})}{\sum \mathrm{U_j}\ \text{from CPLEX}} \quad (13)$$

Figure 9 presents the main effects plot of the various factors on the average % improvement in minimizing the total number of tardy jobs.

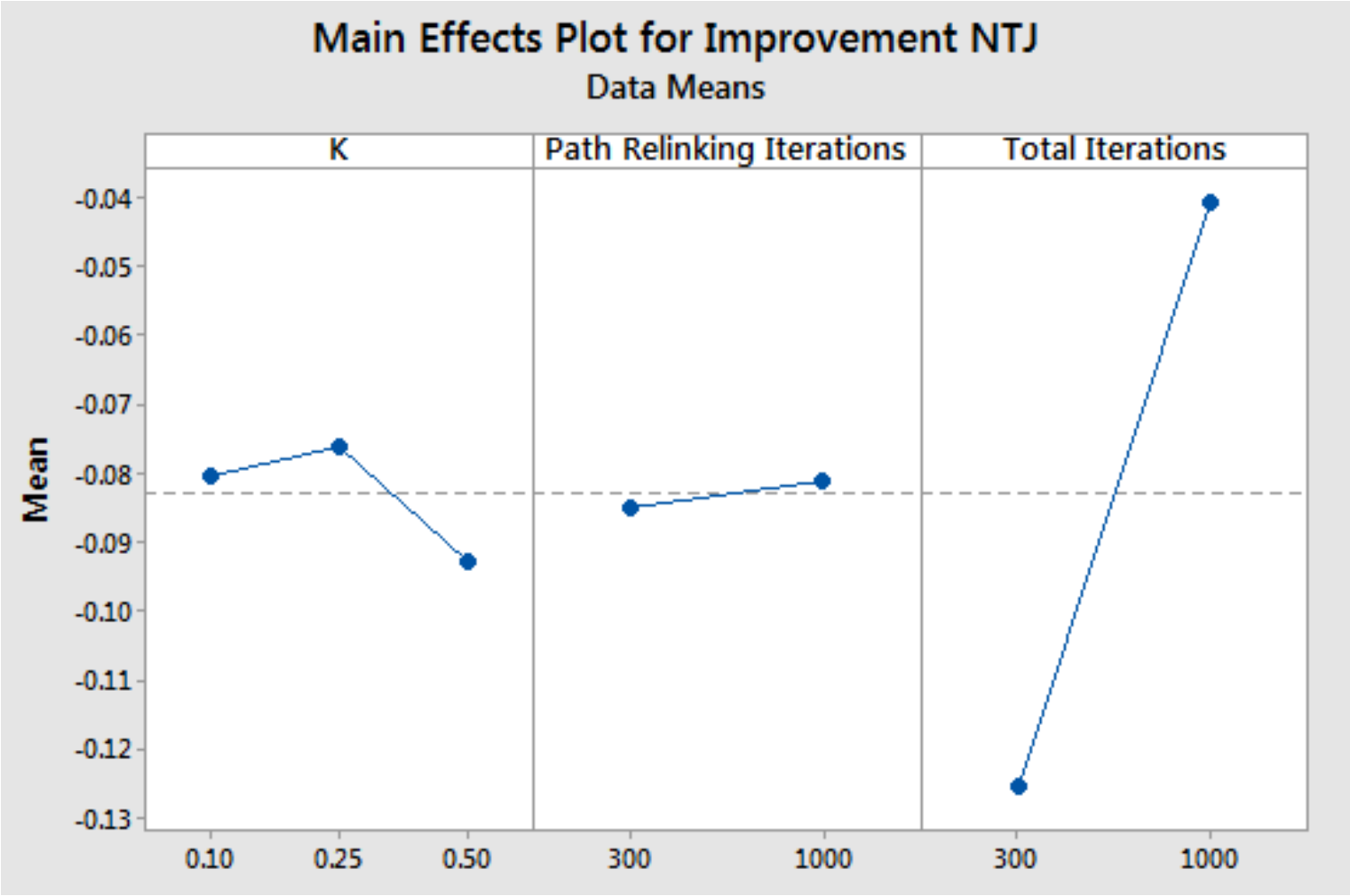

**Figure 9** Main effect plot for GRASP improvement.

The average percentage improvement in minimizing the total number of tardy jobs was higher when the GRASP parameters were set to 1000 iterations, 1000 path relinking iterations, and 0.25 (i.e. 13 for 50 jobs instances) as *k*. As the size of *k* was increased, the improvement in minimizing the total number of tardy jobs slightly decreased. Figure 10 presents the main effects plot of the various factors on the average run time required by GRASP to solve these instances.

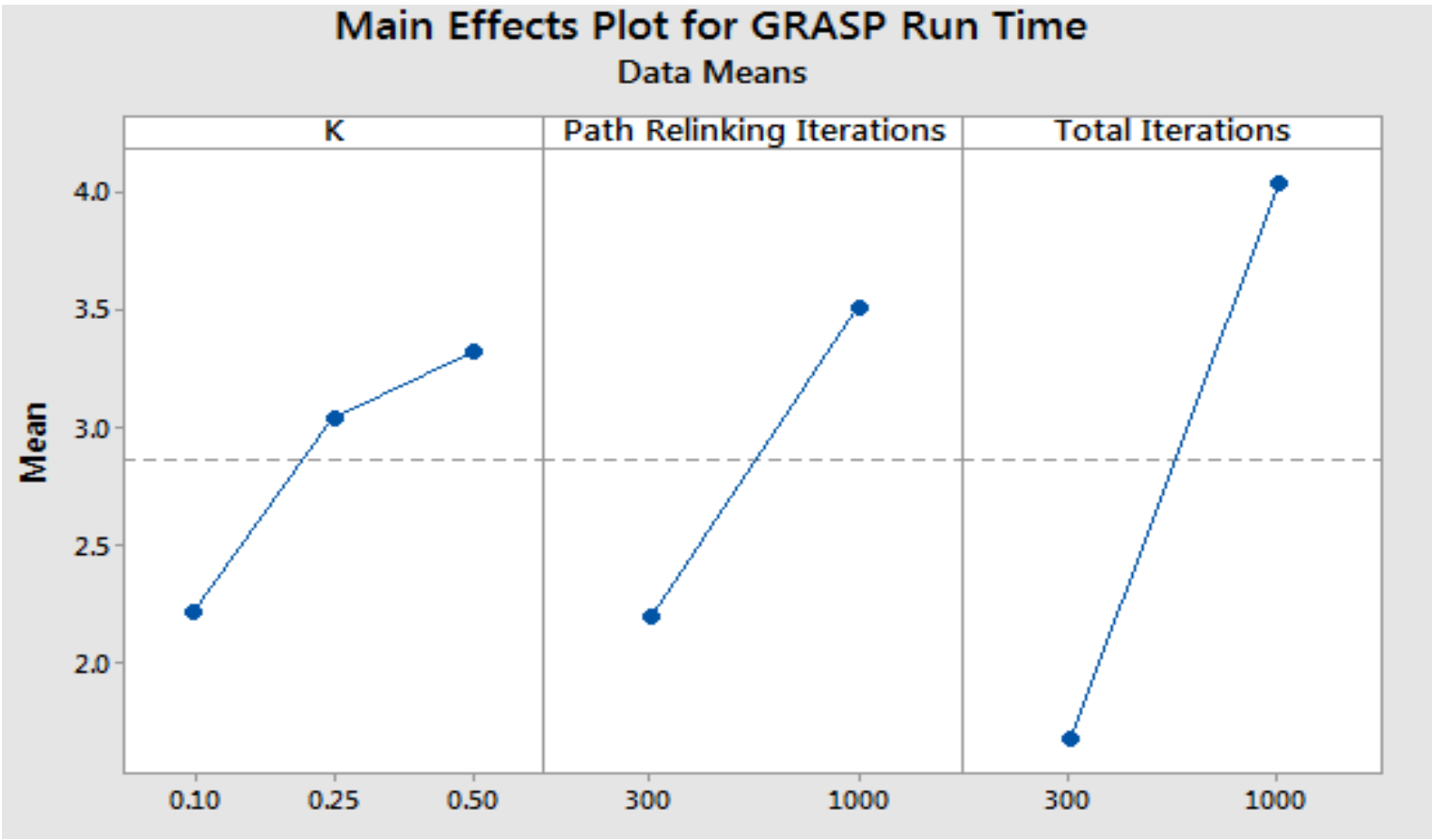

**Figure 10** Main effect plot for GRASP run time.

The run times were shorter when the parameters were set to 300 iterations, 300 path relinking iterations, and 10% of number of jobs (i.e. 5 for 50 job instances) as *k*. Although it requires less time to run GRASP when the total iteration is 300, the solution quality is inferior when compared to the total iteration of 1000. Based on this preliminary research, the GRASP parameters were set as follows: 1000 iterations, 1000 path relinking iterations and *k*=10% of number of jobs. The rest of the problem instances were solved by fixing the GRASP parameters to these values.

The performance of the proposed GRASP algorithm was evaluated by comparing its results with the proposed construction heuristic and the solution from a commercial solver used to solve the formulation proposed in this research. The total number of tardy jobs obtained from the GRASP algorithm was compared to the best results gained from the construction heuristic. The proposed GRASP algorithm and construction heuristic was implemented in MATLAB R2012a.

The mathematical model proposed in this research was solved using the IBM ILOG CPLEX 10.1 solver. Since CPLEX did not converge to an optimum even after running for several hours on most larger problem instances, the best-known solution from CPLEX (after running for about 1800 s) was used for comparison purposes. The total number of tardy jobs from CPLEX was compared with the construction heuristic and GRASP results. Figures 11, 12 and 13 summarize the average percentage improvement in the total number of tardy jobs when GRASP and the construction heuristic are compared with the commercial solver. From Figures 11, 12 and 13, it can be inferred that GRASP reports better solutions than the construction heuristic. Negative percentages indicate that the CPLEX results were better. Positive percentages indicate the CPLEX results were poor compared to either the construction heuristic or the GRASP algorithm. For example, in Figure 11, the GRASP algorithm found a 31.66% average improved solution for the 100 iteration case when compared to the CPLEX solution.

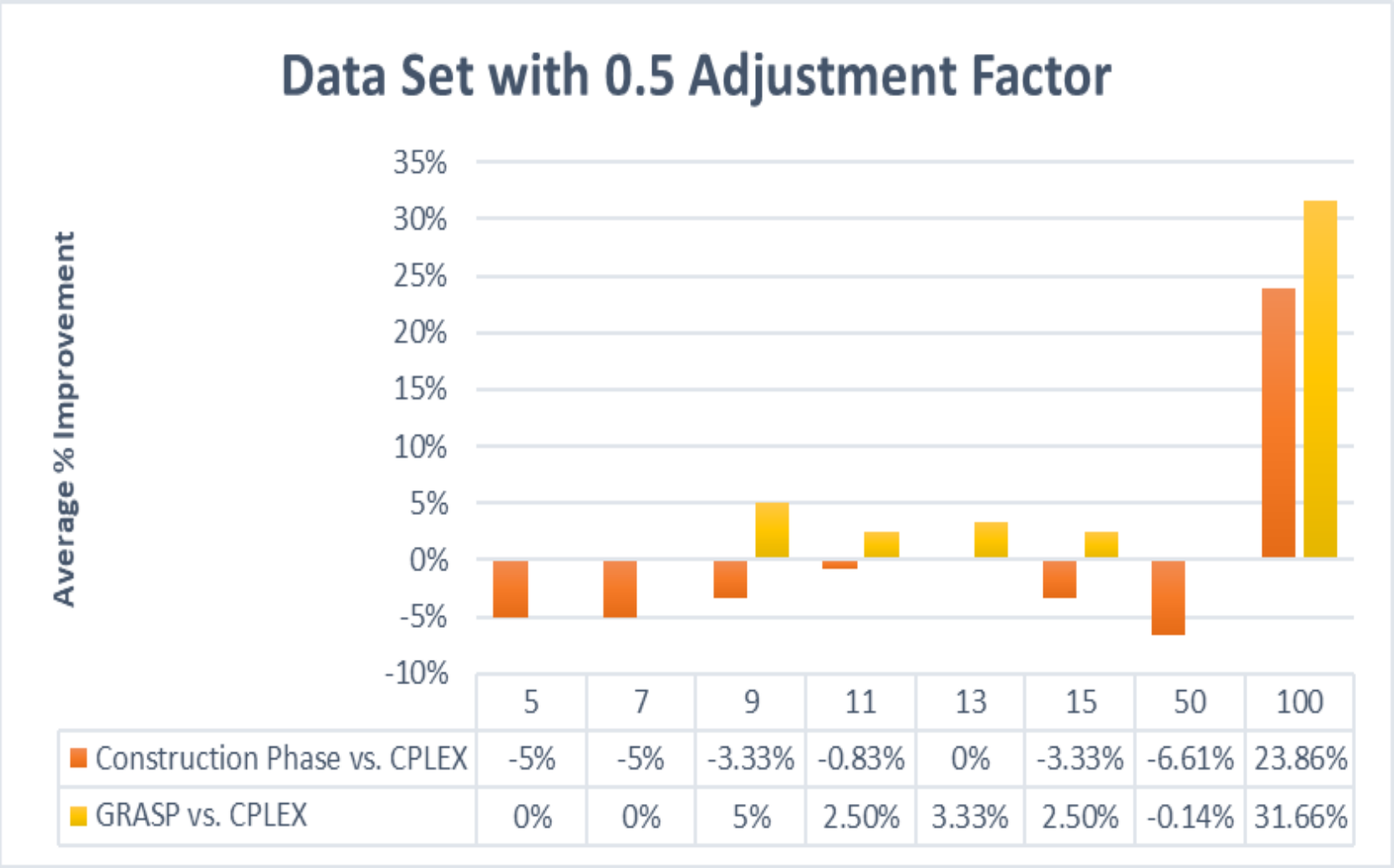

**Figure 11** GRASP vs. construction for data set with 0.5 adjustment factor.

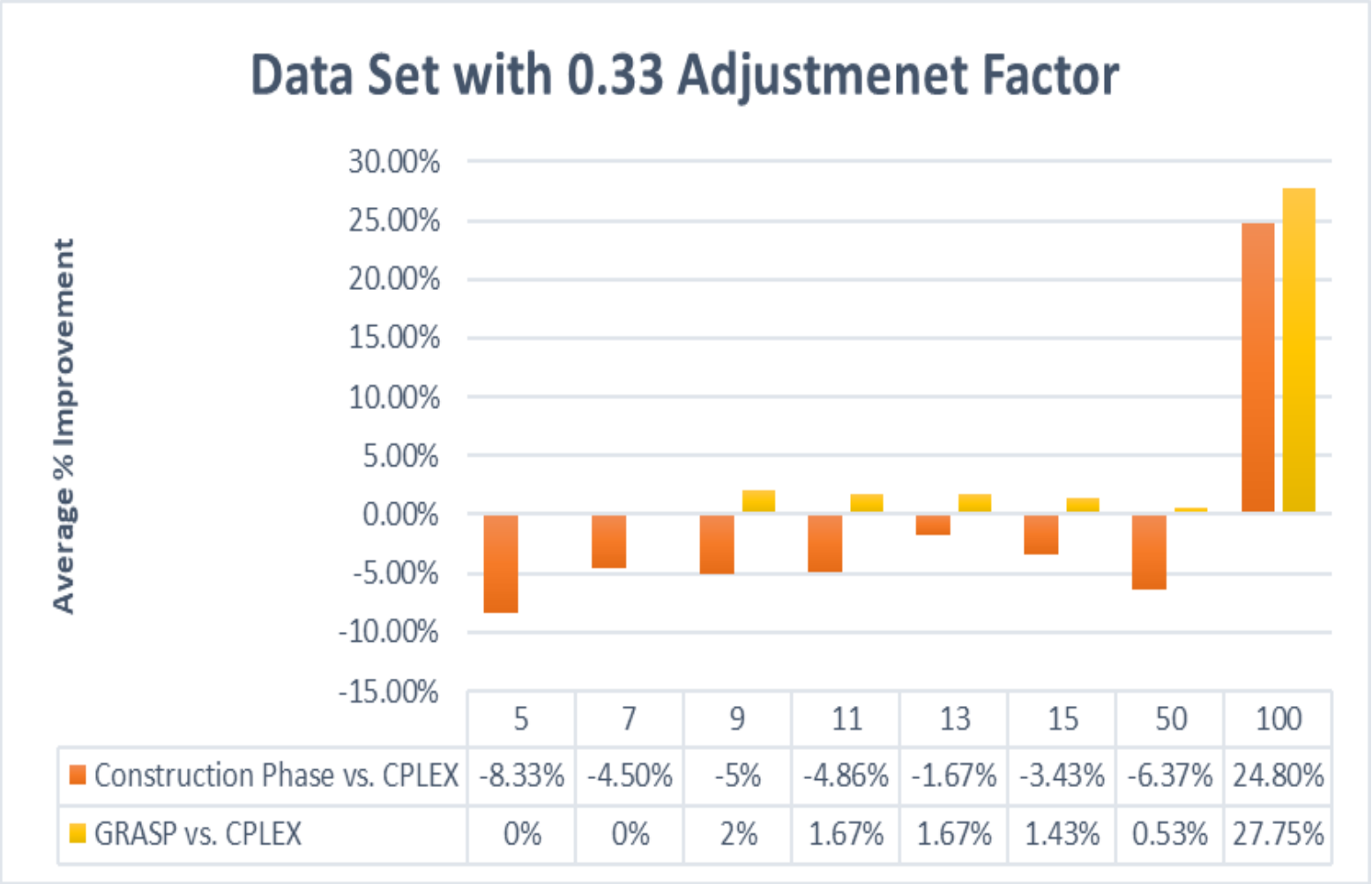

**Figure 12** GRASP vs. construction for data set with 0.33 adjustment factor.

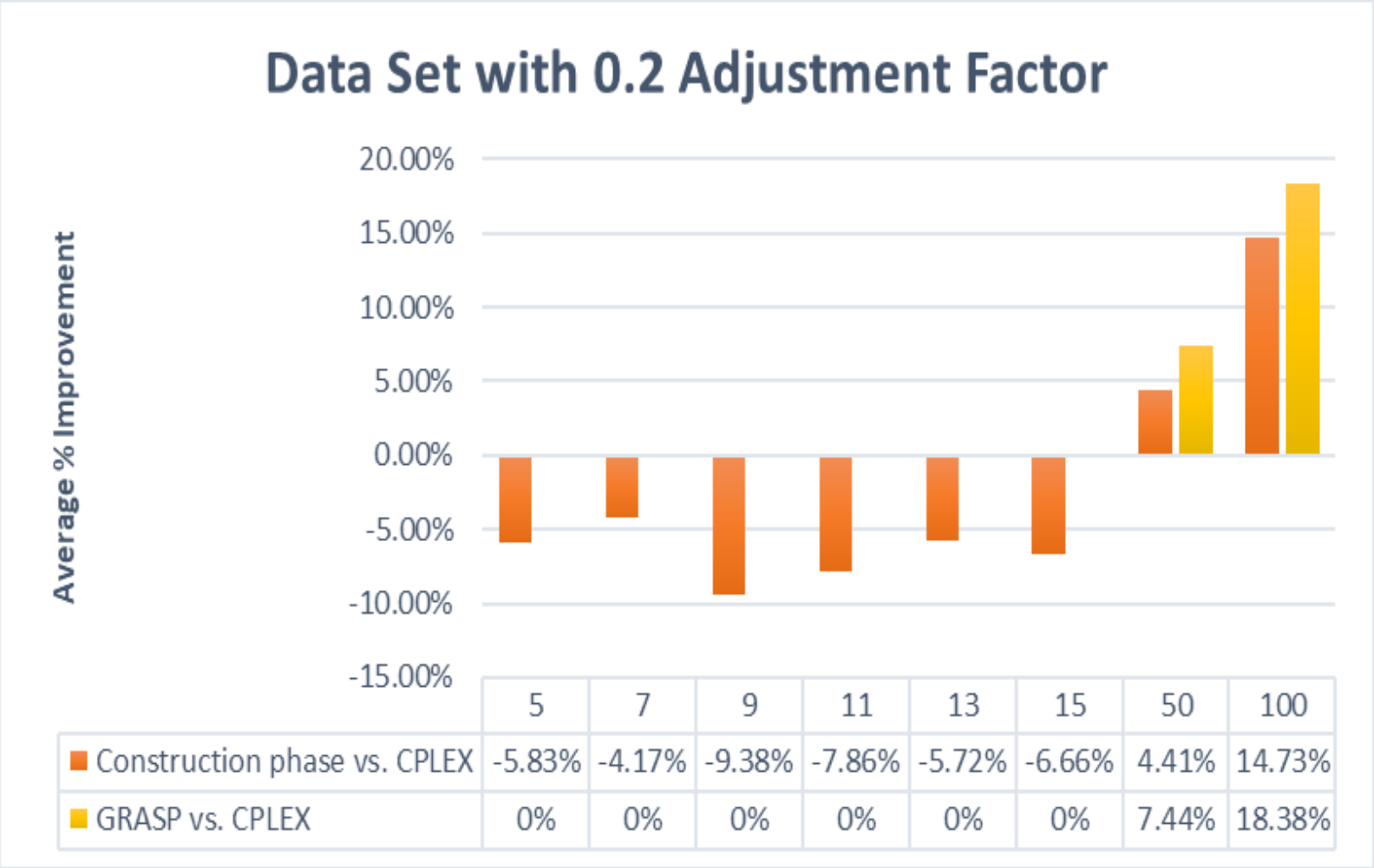

**Figure 13** GRASP vs. construction for data set with 0.2 adjustment factor.

Figures 14, 15 and 16 compare the run times required to solve the problem instances by GRASP and the CPLEX solver. The CPLEX solver was restricted to run for a maximum of 1800s (30 min) as it failed to converge to optimum even after running for several hours. It is evident from these figures that GRASP required shorter run time when compared to CPLEX.

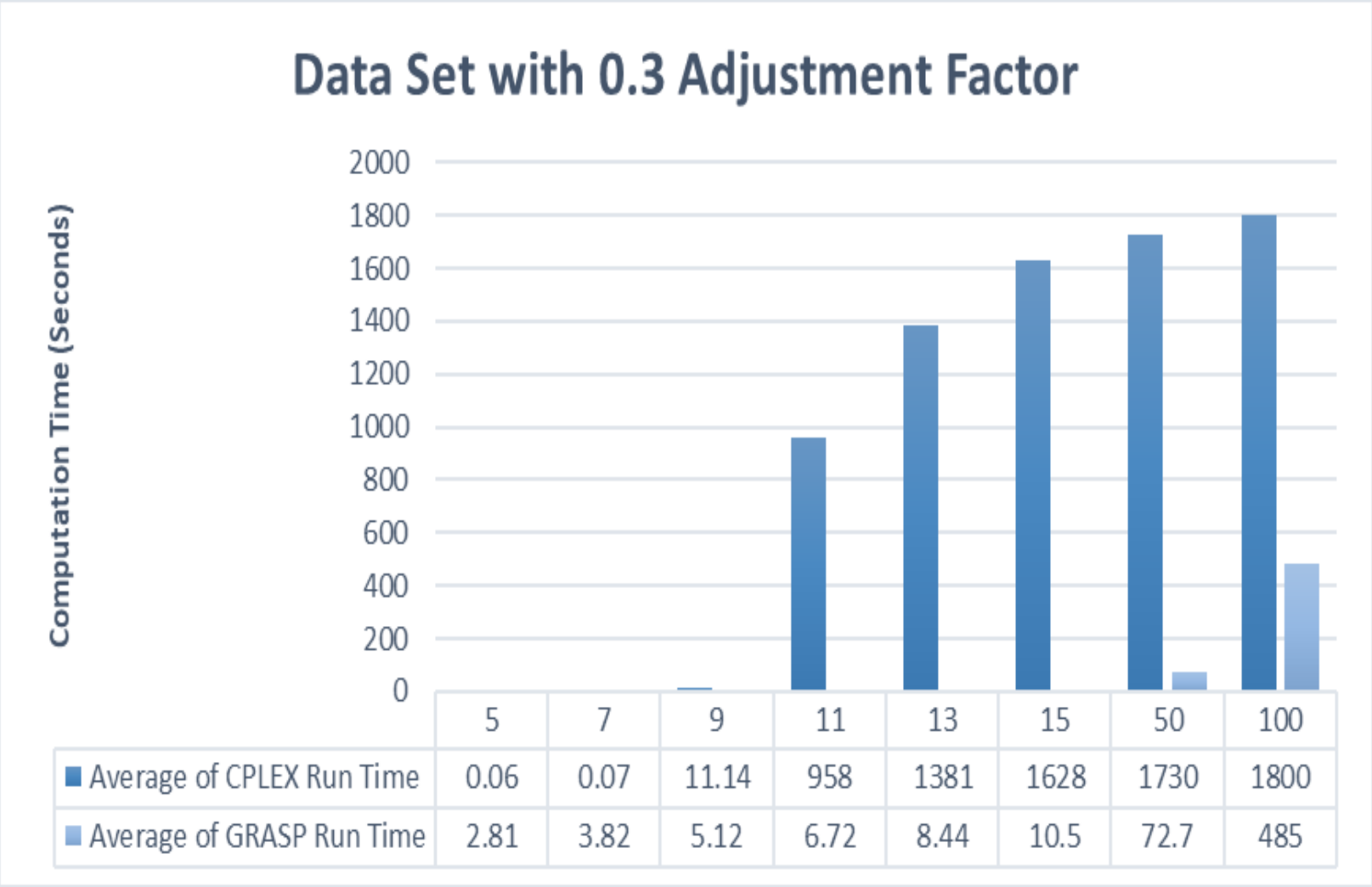

**Figure 14** Average run time comparison CPLEX vs. GRASP for data set with 0.33 adjustment factor.

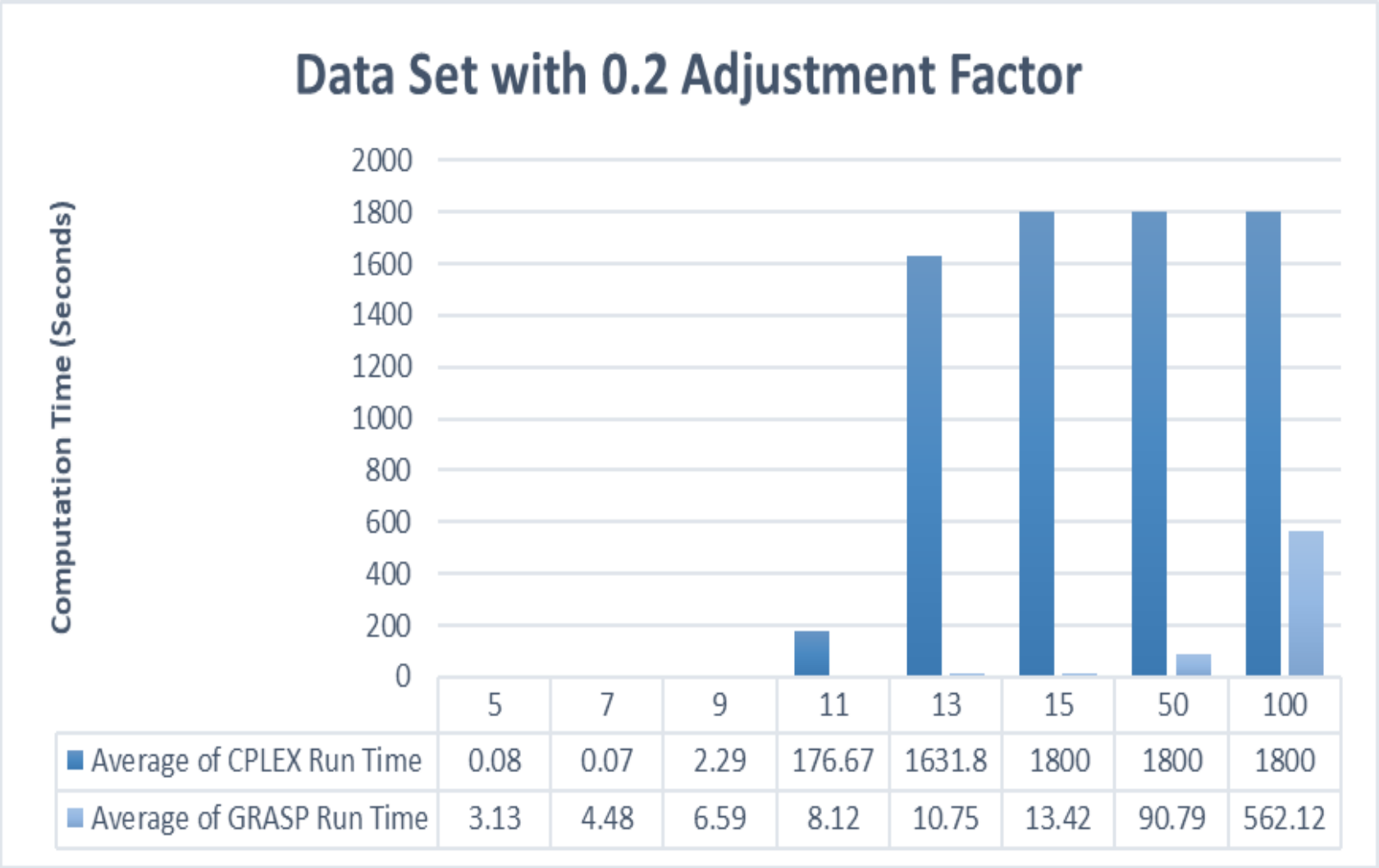

**Figure 15** Average run time comparison CPLEX vs. GRASP for data set with 0.2 adjustment factor.

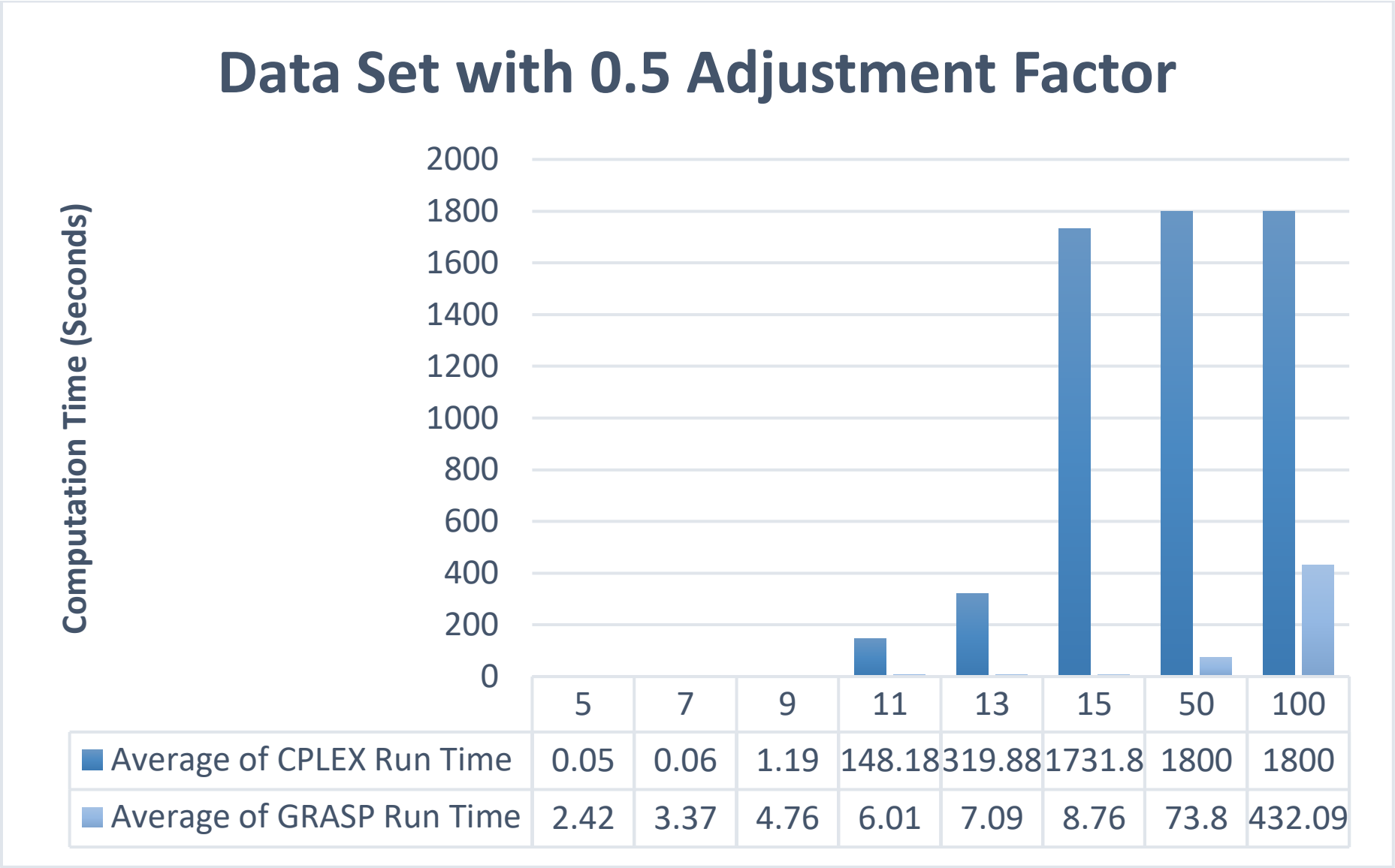

**Figure 16** Average run time comparison CPLEX vs. GRASP for data set with 0.5 adjustment factor.

## 7 Conclusion

This paper presented a GRASP algorithm with path relinking to schedule jobs on a single batch processing machine such that the total number of tardy jobs are minimized. The GRASP algorithm was compared with the proposed construction heuristic and a commercial solver was used to solve a mixed-integer linear program. Through an experimental study, it was shown that GRASP outperformed both construction heuristic and CPLEX in terms of solution quality. GRASP approach also was superior than CPLEX in terms of run time. The research findings indicate that this approach is successful in solving large problem instances and it is practical to implement as it reports good quality solutions in acceptable run time. The contributions of this research would motivate other researchers to consider GRASP for other variants of this and other BPM scheduling problems with a number of tardy jobs objective.